# Simple proofs of Bressoud's and Schur's polynomial versions of the Rogers-Ramanujan identities.

*Johann Cigler*

Fakultät für Mathematik
Universität Wien
A-1090 Wien, Nordbergstraße 15

Johann.Cigler@univie.ac.at

**Abstract**
*We give simple elementary proofs of Bressoud's and Schur's polynomial versions of the Rogers-Ramanujan identities.*

## 1. Bressoud's identity

In [1] George E. Andrews and Kimmo Eriksson gave a simple proof of David Bressoud's ([2]) polynomial version of the Rogers-Ramanujan identities. I want to show that their proof can be further simplified by starting with the identity

$$\sum_{j=-k}^{k} (-1)^j q^{\frac{j(3j-1)}{2}} \begin{bmatrix} n \\ k-j \end{bmatrix}\begin{bmatrix} n \\ k+j \end{bmatrix} = \begin{bmatrix} n \\ k \end{bmatrix}. \tag{1}$$

Ole Warnaar has informed me that this identity has been obtained in [6], Lemma 3.1 as limit case of Rogers' q-Dougall sum. In [6] he already used (1) to prove (a generalization of) Bressoud's identity (11). Christian Krattenthaler has told me that (1) can be considered as limit case of Jackson's q-Dixon summation. It is also a special case of Paule's transformation T1 of [4].

A simple computer proof can be given if we write the left hand side of (1) in the equivalent form

$$f(n,k) = \sum_{j=-k}^{k} (-1)^j q^{\frac{j(3j-1)}{2}} \frac{1+q^j}{2} \begin{bmatrix} n \\ k-j \end{bmatrix}\begin{bmatrix} n \\ k+j \end{bmatrix}.$$

Then the implementation qZeil of the $q$-Zeilberger algorithm gives

$$f(n,k) = \frac{1-q^n}{1-q^{n-k}} f(n-1,k),$$

from which (1) is obvious if we observe that $f(k,k) = 1$.

If you don't trust the computer set $a(n,k,j) = (-1)^j q^{\frac{j(3j-1)}{2}} (1+q^j) \begin{bmatrix} n \\ k-j \end{bmatrix}\begin{bmatrix} n \\ k+j \end{bmatrix}$ and

$b(n,k,j) = \frac{q^{n-k+2j}(1-q^{k-j})(1-q^{n-j-k})}{(1+q^j)(1-q^{n-k})(1-q^n)} a(n,k,j)$ and verify that for $n > k$

$$a(n,k,j) - \frac{1-q^n}{1-q^{n-k}} a(n-1,k,j) = b(n,k,j) - b(n,k,j-1).$$

Here I give an elementary proof of (1) which uses only the recurrence relations for the $q$-binomial coefficients:
To this end let

$$f(n,k) = \sum_{j=-k}^{k} (-1)^j q^{\frac{j(3j-1)}{2}} \begin{bmatrix} n \\ k-j \end{bmatrix} \begin{bmatrix} n \\ k+j \end{bmatrix} = \sum_{j=-k}^{k} (-1)^j q^{\frac{j(3j+1)}{2}} \begin{bmatrix} n \\ k-j \end{bmatrix} \begin{bmatrix} n \\ k+j \end{bmatrix}. \tag{2}$$

From the recurrence relation

$$\begin{bmatrix} n+1 \\ k \end{bmatrix} = \begin{bmatrix} n \\ k \end{bmatrix} + q^{n-k+1} \begin{bmatrix} n \\ k-1 \end{bmatrix} \tag{3}$$

for the $q$-binomial coefficients we also get

$$\sum_j (-1)^j q^{\frac{j(3j-1)}{2}} \begin{bmatrix} n \\ k-j \end{bmatrix} \begin{bmatrix} n+1 \\ k+j \end{bmatrix} = f(n,k). \tag{4}$$

This follows from

$$\sum_j (-1)^j q^{\frac{j(3j-1)}{2}} \begin{bmatrix} n \\ k-j \end{bmatrix} \begin{bmatrix} n+1 \\ k+j \end{bmatrix} = \sum_j (-1)^j q^{\frac{j(3j-1)}{2}} \begin{bmatrix} n \\ k-j \end{bmatrix} \begin{bmatrix} n \\ k+j \end{bmatrix}$$
$$+ \sum_j (-1)^j q^{\frac{j(3j-1)}{2}+n-k-j+1} \begin{bmatrix} n \\ k-j \end{bmatrix} \begin{bmatrix} n \\ k+j-1 \end{bmatrix} = f(n,k) + q^{n-k+1} \sum_j (-1)^j q^{\frac{3j(j-1)}{2}} \begin{bmatrix} n \\ k-j \end{bmatrix} \begin{bmatrix} n \\ k+j-1 \end{bmatrix}.$$

The last sum vanishes, because $j \to -j+1$ defines a sign reversing involution.

The other recurrence relation

$$\begin{bmatrix} n+1 \\ k \end{bmatrix} = q^k \begin{bmatrix} n \\ k \end{bmatrix} + \begin{bmatrix} n \\ k-1 \end{bmatrix} \tag{5}$$

gives

$$\sum_j (-1)^j q^{\frac{j(3j+1)}{2}} \begin{bmatrix} n \\ k-j \end{bmatrix} \begin{bmatrix} n+1 \\ k+j+1 \end{bmatrix} = q^{k+1} \sum_j (-1)^j q^{\frac{3j(j+1)}{2}} \begin{bmatrix} n \\ k-j \end{bmatrix} \begin{bmatrix} n \\ k+j+1 \end{bmatrix} + f(n,k) = f(n,k).$$

Therefore we get

$$f(n,k) = \sum_{j=-k}^{k} (-1)^j q^{\frac{j(3j+1)}{2}} \begin{bmatrix} n \\ k-j \end{bmatrix} \begin{bmatrix} n \\ k+j \end{bmatrix} = \sum_j (-1)^j q^{\frac{j(3j-1)}{2}} \begin{bmatrix} n \\ k-j \end{bmatrix} \begin{bmatrix} n+1 \\ k+j \end{bmatrix}$$
$$= \sum_j (-1)^j q^{\frac{j(3j+1)}{2}} \begin{bmatrix} n \\ k-j \end{bmatrix} \begin{bmatrix} n+1 \\ k+j+1 \end{bmatrix}. \tag{6}$$

This implies

$$f(n+1,k) = \sum_j (-1)^j q^{\frac{j(3j-1)}{2}} \begin{bmatrix} n+1 \\ k-j \end{bmatrix}\begin{bmatrix} n+1 \\ k+j \end{bmatrix} = \sum_j (-1)^j q^{\frac{j(3j-1)}{2}} \begin{bmatrix} n \\ k-j \end{bmatrix}\begin{bmatrix} n+1 \\ k+j \end{bmatrix}$$

$$+ q^{n-k+1} \sum_j (-1)^j q^{\frac{j(3j+1)}{2}} \begin{bmatrix} n \\ k-j-1 \end{bmatrix}\begin{bmatrix} n+1 \\ k+j \end{bmatrix}$$

$$= f(n,k) + q^{n-k+1} \sum_j (-1)^j q^{\frac{j(3j+1)}{2}} \begin{bmatrix} n \\ (k-1)-j \end{bmatrix}\begin{bmatrix} n+1 \\ (k-1)+j+1 \end{bmatrix} = f(n,k) + q^{n-k+1} f(n,k-1).$$

Therefore the sequence $(f(n,k))$ satisfies the recurrence relation (3) for the $q$–binomial coefficients and the corresponding boundary values.

This proves

**Theorem 1**

*The following identities hold:*

$$\sum_{j=-k}^{k} (-1)^j q^{\frac{j(3j-1)}{2}} \begin{bmatrix} n \\ k-j \end{bmatrix}\begin{bmatrix} n \\ k+j \end{bmatrix} = \sum_j (-1)^j q^{\frac{j(3j-1)}{2}} \begin{bmatrix} n \\ k-j \end{bmatrix}\begin{bmatrix} n+1 \\ k+j \end{bmatrix}$$
$$= \sum_j (-1)^j q^{\frac{j(3j+1)}{2}} \begin{bmatrix} n \\ k-j \end{bmatrix}\begin{bmatrix} n+1 \\ k+j+1 \end{bmatrix} = \begin{bmatrix} n \\ k \end{bmatrix}.$$
(7)

From (7) we obtain

$$\sum_k \begin{bmatrix} n \\ k \end{bmatrix} q^{k^2} = \sum_j (-1)^j q^{\frac{j(5j-1)}{2}} \sum_{k \geq |j|} q^{(k-j)(k+j)} \begin{bmatrix} n \\ k-j \end{bmatrix}\begin{bmatrix} n \\ k+j \end{bmatrix}.$$
(8)

The $q$–Vandermonde formula

$$\begin{bmatrix} m+n \\ k \end{bmatrix} = \sum_j \begin{bmatrix} m \\ j \end{bmatrix}\begin{bmatrix} n \\ k-j \end{bmatrix} q^{(m-j)(k-j)}$$
(9)

implies

$$\sum_{k \geq |j|} q^{(k-j)(k+j)} \begin{bmatrix} n \\ k-j \end{bmatrix}\begin{bmatrix} n \\ k+j \end{bmatrix} = \begin{bmatrix} 2n \\ n-2j \end{bmatrix}.$$
(10)

Therefore (8) reduces to Bressoud's identity

$$\sum_k \begin{bmatrix} n \\ k \end{bmatrix} q^{k^2} = \sum_j (-1)^j q^{\frac{j(5j-1)}{2}} \begin{bmatrix} 2n \\ n-2j \end{bmatrix}.$$
(11)

In the same way we get

$$\sum_j(-1)^j q^{\frac{j(3j-3)}{2}}\begin{bmatrix}n\\k-j\end{bmatrix}\begin{bmatrix}n+1\\k+j\end{bmatrix} = q^k\sum_j(-1)^j q^{\frac{j(3j-1)}{2}}\begin{bmatrix}n\\k-j\end{bmatrix}\begin{bmatrix}n\\k+j\end{bmatrix}$$

$$+\sum_j(-1)^j q^{\frac{j(3j-3)}{2}}\begin{bmatrix}n\\k-j\end{bmatrix}\begin{bmatrix}n\\k+j-1\end{bmatrix} = q^k f(n,k) = q^k\begin{bmatrix}n\\k\end{bmatrix}.$$

This implies as above

$$\sum_k\begin{bmatrix}n\\k\end{bmatrix}q^{k^2+k} = \sum_j(-1)^j q^{\frac{j(5j-3)}{2}}\begin{bmatrix}2n+1\\n+1-2j\end{bmatrix}. \tag{12}$$

As is well known (cf. e.g. [1]) by letting $n \to \infty$ in (11) we get the first Rogers-Ramanujan identity

$$\sum_{k\geq 0}\frac{q^{k^2}}{(1-q)(1-q^2)\cdots(1-q^k)} = \frac{1}{\prod_{k\geq 1}(1-q^k)}\sum_{j\in\mathbb{Z}}(-1)^j q^{\frac{j(5j-1)}{2}}. \tag{13}$$

In the same way from (12) we get the second Rogers-Ramanujan identity

$$\sum_{k\geq 0}\frac{q^{k^2+k}}{(1-q)(1-q^2)\cdots(1-q^k)} = \frac{1}{\prod_{k\geq 1}(1-q^k)}\sum_{j\in\mathbb{Z}}(-1)^j q^{\frac{j(5j-3)}{2}}. \tag{14}$$

## 2. Schur's identity

The identity which corresponds to (1) for Schur's polynomial version is

**Theorem 2**

$$g(n,k) = \sum_{j=-k}^{k}(-1)^j q^{\frac{j(3j-1)}{2}}\begin{bmatrix}\left\lfloor\frac{n+j}{2}\right\rfloor\\k-j\end{bmatrix}\begin{bmatrix}\left\lfloor\frac{n-j+1}{2}\right\rfloor\\k+j\end{bmatrix} = \begin{bmatrix}n-k\\k\end{bmatrix}. \tag{15}$$

This identity has been obtained in [3] by other means.

By using (5) we get

$$g(n+2,k) = q^k\sum_{j=-k}^{k}(-1)^j q^{\frac{j(3j-3)}{2}}\begin{bmatrix}\left\lfloor\frac{n+j}{2}\right\rfloor\\k-j\end{bmatrix}\begin{bmatrix}\left\lfloor\frac{n-j+3}{2}\right\rfloor\\k+j\end{bmatrix} + \sum_{j=-k}^{k}(-1)^j q^{\frac{j(3j-1)}{2}}\begin{bmatrix}\left\lfloor\frac{n+j}{2}\right\rfloor\\k-j-1\end{bmatrix}\begin{bmatrix}\left\lfloor\frac{n-j+3}{2}\right\rfloor\\k+j\end{bmatrix}.$$

For the first sum we get again by using (5)

$$\sum_{j=-k}^{k}(-1)^j q^{\frac{j(3j-3)}{2}} \begin{bmatrix} \frac{n+j}{2} \\ k-j \end{bmatrix} \begin{bmatrix} \frac{n-j+3}{2} \\ k+j \end{bmatrix} = q^k \sum_{j=-k}^{k}(-1)^j q^{\frac{j(3j-1)}{2}} \begin{bmatrix} \frac{n+j}{2} \\ k-j \end{bmatrix} \begin{bmatrix} \frac{n-j+1}{2} \\ k+j \end{bmatrix}$$

$$+ \sum_{j=-k}^{k}(-1)^j q^{\frac{j(3j-3)}{2}} \begin{bmatrix} \frac{n+j}{2} \\ k-j \end{bmatrix} \begin{bmatrix} \frac{n-j+1}{2} \\ k+j-1 \end{bmatrix} = q^k g(n,k) + \ell(n,k),$$

where

$$\ell(n,k) = \sum_{j=-k}^{k}(-1)^j q^{\frac{j(3j-3)}{2}} \begin{bmatrix} \frac{n+j}{2} \\ k-j \end{bmatrix} \begin{bmatrix} \frac{n-j+1}{2} \\ k+j-1 \end{bmatrix} = 0,$$

because $j \to -j+1$ induces a sign reversing involution.

Therefore we have

$$\sum_{j=-k}^{k}(-1)^j q^{\frac{j(3j-3)}{2}} \begin{bmatrix} \frac{n+j}{2} \\ k-j \end{bmatrix} \begin{bmatrix} \frac{n-j+3}{2} \\ k+j \end{bmatrix} = q^k g(n,k). \tag{16}$$

The second term in the above formula gives

$$\sum_{j=-k}^{k}(-1)^j q^{\frac{j(3j-1)}{2}} \begin{bmatrix} \frac{n+j}{2} \\ k-j-1 \end{bmatrix} \begin{bmatrix} \frac{n-j+3}{2} \\ k+j \end{bmatrix} = q^k \sum_{j=-k}^{k}(-1)^j q^{\frac{j(3j+1)}{2}} \begin{bmatrix} \frac{n+j}{2} \\ k-j-1 \end{bmatrix} \begin{bmatrix} \frac{n-j+1}{2} \\ k+j \end{bmatrix}$$

$$+ \sum_{j=-k}^{k}(-1)^j q^{\frac{j(3j-1)}{2}} \begin{bmatrix} \frac{n+j}{2} \\ k-j-1 \end{bmatrix} \begin{bmatrix} \frac{n-j+1}{2} \\ k+j-1 \end{bmatrix} = q^k \sum_{j=-k}^{k}(-1)^j q^{\frac{j(3j-1)}{2}} \begin{bmatrix} \frac{n+j+1}{2} \\ k-j \end{bmatrix} \begin{bmatrix} \frac{n-j}{2} \\ k+j-1 \end{bmatrix} + g(n,k-1).$$

Let now

$$h(n,k) = \sum_{j=-k}^{k}(-1)^j q^{\frac{j(3j-1)}{2}} \begin{bmatrix} \frac{n+j+1}{2} \\ k-j \end{bmatrix} \begin{bmatrix} \frac{n-j}{2} \\ k+j-1 \end{bmatrix}.$$

Then

$$h(n,k) = \sum_{j=-k}^{k} (-1)^j q^{\frac{j(3j-1)}{2}} \left[\begin{array}{c} \left\lfloor\frac{n+j+1}{2}\right\rfloor \\ k-j \end{array}\right] \left[\begin{array}{c} \left\lfloor\frac{n-j}{2}\right\rfloor \\ k+j-1 \end{array}\right]$$

$$= q^k \sum_{j=-k}^{k} (-1)^j q^{\frac{j(3j-3)}{2}} \left[\begin{array}{c} \left\lfloor\frac{n+j-1}{2}\right\rfloor \\ k-j \end{array}\right] \left[\begin{array}{c} \left\lfloor\frac{n-j}{2}\right\rfloor \\ k+j-1 \end{array}\right] + \sum_{j=-k}^{k} (-1)^j q^{\frac{j(3j-1)}{2}} \left[\begin{array}{c} \left\lfloor\frac{n+j-1}{2}\right\rfloor \\ k-j-1 \end{array}\right] \left[\begin{array}{c} \left\lfloor\frac{n-j}{2}\right\rfloor \\ k+j-1 \end{array}\right]$$

$$= q^k \ell(n-1,k) + g(n-1,k-1) = g(n-1,k-1).$$

Therefore we get the recursion

$$g(n+2,k) = q^{2k} g(n,k) + q^k g(n-1,k-1) + g(n,k-1). \qquad (17)$$

It is easy to verify that $g(n,0) = 1 = \left[\begin{array}{c} n-0 \\ 0 \end{array}\right]$, $g(k,k) = 0 = \left[\begin{array}{c} k-k \\ k \end{array}\right]$ for $k \geq 1$ and

$$g(k+1,k) = \left[\begin{array}{c} k+1-k \\ k \end{array}\right] = 0 \text{ for } k \geq 2.$$

By this recurrence and the initial values $g(n,k)$ is uniquely determined for all $n \geq k$.

Since

$$\left[\begin{array}{c} n+2-k \\ k \end{array}\right] = q^k \left[\begin{array}{c} n+1-k \\ k \end{array}\right] + \left[\begin{array}{c} n+1-k \\ k-1 \end{array}\right] = q^{2k} \left[\begin{array}{c} n-k \\ k \end{array}\right] + q^k \left[\begin{array}{c} n-1-(k-1) \\ k-1 \end{array}\right] + \left[\begin{array}{c} n-(k-1) \\ k-1 \end{array}\right]$$

we see that $g(n,k) = \left[\begin{array}{c} n-k \\ k \end{array}\right]$ for all $n \geq k$.

By summing over all $k$ and using the $q-$Vandermonde formula we get

$$\sum_{k=0}^{n} q^{k^2} \left[\begin{array}{c} n-k \\ k \end{array}\right] = \sum_{k=0}^{n} q^{k^2} \sum_{j=-k}^{k} (-1)^j q^{\frac{j(3j-1)}{2}} \left[\begin{array}{c} \left\lfloor\frac{n+j}{2}\right\rfloor \\ k-j \end{array}\right] \left[\begin{array}{c} \left\lfloor\frac{n-j+1}{2}\right\rfloor \\ k+j \end{array}\right]$$

$$= \sum_{j=-n}^{n} (-1)^j q^{\frac{j(5j-1)}{2}} \sum_{k \geq |j|} q^{(k-j)(k+j)} \left[\begin{array}{c} \left\lfloor\frac{n+j}{2}\right\rfloor \\ k-j \end{array}\right] \left[\begin{array}{c} \left\lfloor\frac{n-j+1}{2}\right\rfloor \\ k+j \end{array}\right] = \sum_{j=-n}^{n} (-1)^j q^{\frac{j(5j-1)}{2}} \left[\begin{array}{c} n \\ \left\lfloor\frac{n+5j}{2}\right\rfloor \end{array}\right].$$

This gives Schur's ([5]) polynomial version of the first Rogers-Ramanujan identity

$$\sum_{k=0}^{n} q^{k^2} \left[\begin{array}{c} n-k \\ k \end{array}\right] = \sum_{j=-n}^{n} (-1)^j q^{\frac{j(5j-1)}{2}} \left[\begin{array}{c} n \\ \left\lfloor\frac{n+5j}{2}\right\rfloor \end{array}\right]. \qquad (18)$$

In the same way from (16) we get

$$\sum_{k=0}^{n} q^{k^2+k} \begin{bmatrix} n-k \\ k \end{bmatrix} = \sum_{k=0}^{n} q^{k^2} \sum_{j=-k}^{k} (-1)^j q^{\frac{j(3j-3)}{2}} \begin{bmatrix} \left\lfloor \frac{n+j}{2} \right\rfloor \\ k-j \end{bmatrix} \begin{bmatrix} \left\lfloor \frac{n-j+3}{2} \right\rfloor \\ k+j \end{bmatrix}$$

$$= \sum_{j=-n}^{n} (-1)^j q^{\frac{j(5j-3)}{2}} \sum_{k \geq |j|} q^{(k-j)(k+j)} \begin{bmatrix} \left\lfloor \frac{n+j}{2} \right\rfloor \\ k-j \end{bmatrix} \begin{bmatrix} \left\lfloor \frac{n-j+3}{2} \right\rfloor \\ k+j \end{bmatrix} = \sum_{j=-n}^{n} (-1)^j q^{\frac{j(5j-3)}{2}} \begin{bmatrix} n+1 \\ \left\lfloor \frac{n-5j+3}{2} \right\rfloor \end{bmatrix}.$$

This is Schur's polynomial version of the second Rogers-Ramanujan identity, which is usually written in the form

$$\sum_{k=0}^{n-1} q^{k^2+k} \begin{bmatrix} n-k-1 \\ k \end{bmatrix} = \sum_{j=-n}^{n} (-1)^j q^{\frac{j(5j-3)}{2}} \begin{bmatrix} n \\ \left\lfloor \frac{n-5j+2}{2} \right\rfloor \end{bmatrix}. \tag{19}$$